%mvnm.tex: a Plain TeX file by Shalosh B. Ekhad and Doron Zeilberger
%A Linear Time, and Constant Space, algorithm to Compute the Mixed Moments of the Multivariate Normal Distributions
%begin macros

\baselineskip=14pt
\parskip=10pt

\magnification=\magstephalf

\def\1{{\overline{1}}}
\def\2{{\overline{2}}}
\parindent=0pt
\overfullrule=0in

\def\frac#1#2{{#1 \over #2}}
%\headline={\rm  \ifodd\pageno  \RightHead  \else  \LeftHead  \fi}
%\def\RightHead{\centerline{
%Title
%}}
%\def\LeftHead{ \centerline{Doron Zeilberger}}
%end macros
\centerline
{\bf 
A Linear Time and Constant Space Algorithm to Compute the Mixed Moments 
}
\centerline
{\bf of the Multivariate Normal Distributions}
\bigskip
\centerline
{\it Shalosh B. EKHAD and Doron ZEILBERGER}
\bigskip

{\bf Abstract}:
Using recurrences gotten from the Apagodu-Zeilberger Multivariate Almkvist-Zeilberger algorithm, we present a linear-time and constant-space algorithm to compute 
the general mixed moments of the $k$-variate general normal distribution, with any covariance matrix, for any specific $k$. 
Besides their obvious importance in statistics, they are also very significant in enumerative combinatorics, since, when the entries of the covariance
matrix remain symbolic, they enable us to count in how
many ways, in a species with $k$ different genders, a bunch of individuals can all get married, keeping track of the different kinds of the $k(k-1)/2$ possible
heterosexual marriages, and the $k$ possible same-sex marriages.
We completely implement our algorithm (with an accompanying Maple package, MVNM.txt) 
for the bivariate and trivariate cases (and hence taking care of our own $2$-sex society and a putative $3$-sex society), but alas, the actual recurrences for larger k took too long for us to compute. We leave them as computational challenges.

{\bf Maple Package}

This article is accompanied by a Maple package, {\tt MVNM.txt}, available from

{\tt https://sites.math.rutgers.edu/\~{}zeilberg/tokhniot/MVNM.txt} \quad .

The web-page of this article,

{\tt https://sites.math.rutgers.edu/\~{}zeilberg/mamarim/mamarimhtml/mvnm.html} \quad ,

contains input and output files, referred to in this paper.

{\bf The multivariate Normal Distribution}

Recall that the {\it probability density function} (see [T] and [Wik]) of the {\it multivariate normal distribution} with mean ${\bf 0}$ and (symmetric) covariance matrix
${\bf C}= (c_{ij})_{1\leq i,j \leq k}$ is
$$
f_{\bf C} ({\bf x}) :=
\frac{ 
e^{-\frac{1}{2} {\bf x}^T {\bf C}^{-1} {\bf x}}}
{\sqrt{(2 \pi)^k \det {\bf C}}} \quad .
$$

By simple rescaling we can always assume that all the variances are $1$, in other words, that the entries of the main diagonal of ${\bf C}$ are all $1$.

We are interested in fast computation of the {\bf mixed moments}
$$
M_{{\bf C}} (m_1, \cdots , m_k) := \int_{R^k} x_1^{m_1} \cdots x_k^{m_k} \, f_{\bf C} (x_1, \dots, x_k) \, dx_1 \cdots dx_k \quad .
$$

One way (not a good one!) to compute these moments, for any specific $(m_1, \dots, m_k)$ is to diagonalize ${\bf C}$, make a change of variables
and compute an integral of the form
$$
\int_{R^k}  \prod_{i=1}^k \left ( \sum_{j=0}^{k} b_{ij}x_j \right)^{m_i} e^{-\frac{1}{2}(x_1^2+ ... +x_k^2)} \, dx_1 \cdots dx_k \quad.
$$
Then expand  $\prod_{i=1}^k \left ( \sum_{j=0}^{k} b_{ij}x_j \right)^{m_i}$ and use the fact that $\int_{-\infty}^{\infty} e^{-x^2/2} x^{r}\, dx $ is $0$ if $r$ is odd and 
$\sqrt{2\pi}\, \cdot \, \frac{r!}{2^{r/2} (r/2)!)}$ if $r$ is even.

A much better way is via the {\bf moment generating function} ([Wik][T])
$$
\sum_{0 \leq m_1, \dots, m_k < \infty} M_{\bf C}(m_1, \dots, m_k) \frac{t_1^{m_1} \cdots t_k^{m_k}}{m_1! \cdots m_k!} \, = \,
e^{\frac{1}{2} (\sum_{1 \leq i,j \leq k} t_i c_{ij} t_j)} \quad .
$$

This is implemented in procedure {\tt MOMd} in the Maple package {\tt MVNM.txt} mentioned above.
For example to get the $(3,3,3,3)$-mixed moment for the {\bf generic} four-variate normal distribution, with
a general (symbolic) covariance matrix
$$
\left ( \matrix{1    &  c12  &  c13 & c14 \cr
                c12  &  1    &  c23 & c24 \cr
                c13  & c23   &   1  & c34 \cr
                c14  & c24   &  c34 & 1
               }
\right ) \quad ,
$$
enter

{\tt lprint(MOMd([[1,c12,c13,c14],[c12,1,c23,c24],[c13,c23,1,c34],[c14,c24,c34,1]],[3,3,3,3]));}

Defining $M:=max(m_1, \dots, m_k)$,  this requires  $O(M^k)$ time and memory.

Another way is to to use the fact that

$$
\int_{R^k} 
\frac{\partial}{\partial x_1} 
\left ( x_1^{m_1} \cdots x_k^{m_k} \, f_{\bf C} (x_1, \dots, x_k) \right ) \, dx_1 \cdots dx_k \, = \, 0 \quad .
$$

Using the product and the chain rule, and expanding, one gets a certain {\it mixed} recurrence, requiring to compute all the (up to) $m_1 \cdots m_k$ `previous' values, 
requiring, again $O(M)^k)$ memory and time.

But thanks to the {\bf Apagodu-Zeilberger} [ApZ] multivariate extension of the {\bf Almkvist-Zeilberger} [AlZ]
algorithm there exist {\bf pure} recurrences, with polynomial coefficients in $m_1, \dots, m_k$,
in {\bf each} of the discrete coordinate directions. The ones for $k=2$ are fairly simple (they are essentially second-order), but the ones for $k=3$ are already very
complicated. But {\it once found} (and we did find them!) this enables a linear-time and constant-space algorithm for computing any $(m_1,m_2,m_3)$-mixed moment.
The recurrences are too complicated to be typeset here, but can be read from the Maple source-code of procedure {\tt MOM3} in our Maple package.

The syntax is

{\tt MOM3(c12,c13,c23,[m1,m2,m3]);} \quad .

For example to get the $(10,10,10)$ mixed moment as a polynomial in the  symbols $c12,c13,c23$, type:

{\tt MOM3(c12,c13,c23,[10,10,10]);} \quad .

This should (and does!) give the same answer as

{\tt MOMd([[1,c12,c13],[c12,1,c23],[c13,c23,1]],[10,10,10]);} \quad .

To really appreciate the superiority of our algorithm, using {\tt MOM3}, over the straightforward {\tt MOM3d}, try, for example

{\tt restart: read `MVNM.txt`: t0:=time():lu1:=MOM3(1/2,1/3,1/4,[570,560,750]); time()-t0;} \quad ,

that would give you the very complicated {\tt lu1} in 2.56 seconds, while

{\tt t0:=time(): lu2:=MOMd([[1,1/2,1/3],[1/2,1,1/4],[1/3,1/4,1]],[570,560,750]);} \quad ,

would confirm that {\tt lu2} and {\tt lu1} are the same (good check!), but it takes 631.007 seconds.

{\bf Warning:} Don't even try to use floating-points! You would get garbage, due to the complexity of the calculations that accumulate the round-off errors.
Both ways would give you erroneous answers unless {\tt Digits} is set very high.

If you keep $c12,c12,c23$ symbolic, the superiority of {\tt MOM3} over {\tt MOMd} is even more apparent.

{\tt restart: read `MVNM.txt`: time(MOM3(c12,c13,c23,[100,50,40]));} \quad ,
is less than $12$ seconds, while doing the same things with {\tt MOMd} takes $100$ times longer!

{\bf Why this is also Important in Enumerative Combinatorics?}

Using what Herb Wilf [Wil] used to call {\it generatiningfunctionlogy} it is easy to see that, when the entries of the covariance matrix $C$ are kept {\it symbolic},
then for the bivariate case, the coefficient of $c^r$ in $M_{[[1,c],[c,1]]}(m1,m2)$ is the exact number of ways that $m1$ men and $m2$ women can all get married
and there are exactly $r$ heterosexual marriages. The coefficiet of
$$
c_{12}^{a_{12}}\,  c_{13}^{a_{23}} \,  c_{23}^{a_{23}}  \quad,
$$
in the polynomial
$$
M_{[[1,c_{12},c_{13}],[c_{12},1,c_{23}],[c_{13},c_{23},1]]}(m_1,m_2,m_3) \quad,
$$
is the exact number of ways that, in a $3$-gender society, with genders $S_1$,  $S_2$, $S_3$, 
that $m_1$ individuals of gender $S_1$,   $m_2$ individuals of gender $S_2$, and  $m_3$ individuals of gender $S_3$,  can {\bf all} get married
(note that you need their total number, $m_1+m_2+m_3$ to be even, or else it is not possible) where there were exactly $a_{12}$  $\{S1,S2\}$ marriages, 
$a_{13}$ $\{S1,S3\}$ marriages,  and $a_{23}$ $\{S2,S3\}$ marriages.

For example, if you want to know the number of ways  $300$ men and $200$ women can get married where there were exactly $100$ heterosexual weddings (and hence $150$ same-sex marriages), type:

{\tt coeff(MOM2(c,[300,200]),c,100);} \quad ,

to get a certain $564$-digit integer.

If you want to know, in a $3$-gender society, the {\bf exact} number of ways that   $20$ individuals of gender S1, $20$ individuals of gender S2, and $20$ individuals of gender S3, can get married
(so altogether there are $30$ weddings) with $9$ $\{S1,S2\}$ weddings, $7$  $\{S1,S3\}$ weddings, and $5$ $\{S2,S3\}$ weddings (and hence $30-9-7-5=9$ same-sex marriages), type:

{\tt coeff(coeff(coeff(MOM3(c12,c13,c23,[20,20,20]),c12,9),c13,7),c23,5);} \quad,

getting, in {\it 0.533} seconds, that the number is:
$$
444975998773143505634352562176000000000 \quad .
$$

{\bf Sample Data}

To see the list of lists of lists of polynomials in {\tt c12,c13,c23}, let's call  it $L$, such that 

{\tt L[m1][m2][m3]}

is the {\tt (m1,m2,m3)}-mixed moment of the trivariate normal distribution with covariance matrix {\tt [[1,c12,c13],[c12,1,c13],[c13,c23,1]]} for $1 \leq m1,m2,m3 \leq 20$
look at the output file

{\tt https://sites.math.rutgers.edu/\~{}zeilberg/tokhniot/oMVNM1.txt}  \quad.

To see the first $35$ diagonal mixed moments (i.e. up to the $(70,70,70)$ mixed moment), see

{\tt https://sites.math.rutgers.edu/\~{}zeilberg/tokhniot/oMVNM2.txt}  \quad.

Enjoy!

The recurrences for four dimensions took too long for us, and we leave them as computational challenges. Perhaps they can be done
with Christoph Koutschan's [K] very powerful Mathematica package? 

{\bf References}

[AlZ]  Gert Almkvist and Doron Zeilberger, {\it The method of differentiating Under The integral sign}, J. Symbolic Computation {\bf 10} (1990), 571-591. \hfill \break
{\tt https://sites.math.rutgers.edu/\~{}zeilberg/mamarim/mamarimPDF/duis.pdf} \quad .

[ApZ] Moa Apagodu and Doron Zeilberger,
{\it Multi-Variable Zeilberger and Almkvist-Zeilberger Algorithms and the Sharpening of Wilf-Zeilberger Theory},
Adv. Appl. Math. {\bf 37} (2006), 139-152. [Special issue in honor of Amitai Regev] \hfill\break
{\tt https://sites.math.rutgers.edu/\~{}zeilberg/mamarim/mamarimhtml/multiZ.html} \quad .

[K] Christoph Koutschan, {\it   Advanced applications of the holonomic systems approach}, 
PhD thesis, Research Institute for Symbolic Computation (RISC), Johannes Kepler University, Linz, Austria, 2009.\hfill\break
{\tt http://www.koutschan.de/publ/Koutschan09/thesisKoutschan.pdf}, \hfill\break
{\tt http://www.risc.jku.at/research/combinat/software/HolonomicFunctions/} \quad.

[T] Y.L. Tong, {\it ``The multivariate normal distribution. Springer Series in Statistics''}. New York: Springer-Verlag, 1990.

[Wik] Wikipedia contributors. {\it Multivariate normal distribution},  Wikipedia, The Free Encyclopedia. Wikipedia, The Free Encyclopedia, 5 Feb. 2022. Web. 8 Feb. 2022.

[Wil] Herbert S. Wilf, {\it ``generatingfunctionology},  Academic Press, 1990. Second Edition: 1994; Third edition : 2005 (CRC Press). Freely downloadable from: \hfill \break
{\tt https://www2.math.upenn.edu/\~{}wilf/gfologyLinked2.pdf}

\bigskip
\hrule
\bigskip
Shalosh B. Ekhad and Doron Zeilberger, Department of Mathematics, Rutgers University (New Brunswick), Hill Center-Busch Campus, 110 Frelinghuysen
Rd., Piscataway, NJ 08854-8019, USA. \hfill\break
Email: {\tt [ShaloshBEkhad, DoronZeil] at gmail dot com}   \quad .

{\bf Exclusively published in the Personal Journal of Shalosh B. Ekhad and Doron Zeilberger and arxiv.org}

{\bf Feb. 20, 2022}

\end